\documentstyle[11pt]{article}
\newtheorem{prop}{Proposition}[section]
\newtheorem{th}[prop]{Theorem}

\newtheorem{lm}[prop]{Lemma}

\newtheorem{remark}[prop]{Remark}

\newcommand{\Ker}{{\rm{Ker}}}
\newcommand{\bsquare}{\hbox{\rule{6pt}{6pt}}}
\newcommand{\proof}[1]{\noindent{\bf Proof}\hspace{0.3cm}{#1}\hfill\bsquare 
\vspace{0.5cm}\par}

\date{}

\begin{document}

\title{Rational curves on the supersingular K3 surface with Artin invariant 1 in characteristic 3}

\author{Toshiyuki Katsura\thanks{Partially supported by JSPS Grant-in-Aid (S), No 19104001} and Shigeyuki Kond\=o\thanks{Partially supported by JSPS Grant-in-Aid (S), No 22224001}}

\maketitle

\begin{abstract}
We show the existence of $112$ non-singular rational curves on the supersingular K3 surface with Artin invariant 1 in characteristic 3
by several ways. Using these rational curves, we have a $(16)_{10}$-configuration 
and a $(280_{4}, 112_{10})$-configuration
on the K3 surface. Moreover we study the Picard lattice by using the theory of the Leech lattice.   The $112$ non-singular rational curves 
correspond to $112$ Leech roots.

\end{abstract}


\section{Introduction}

A supersingular K3 surface (in the sense of Shioda \cite{Sh}) is a K3 surface with the Picard group of rank 22.
This occurs only if the ground field $k$ has a positive characteristic $p$.  Let $S_Y$ be the Picard lattice of
a supersingular K3 surface $Y$ and $S_Y^*$ the dual of $S_Y$.  By a result of Artin \cite{A}, 
$S_Y^*/S_Y$ is a $p$-elementary abelian group $({\bf Z}/p{\bf Z})^{2\sigma}$.  The $\sigma$ is called Artin invariant.
Supersingular K3 surfaces in characteristic $p$ with Artin invariant $\sigma$ form a family of dimension $\sigma -1$.
Moreover Ogus \cite{O1}, \cite{O2} proved the uniqueness of supersingular K3 surface with Artin invariant 1 up to isomorphisms (see also Rudakov and Shafarevich \cite{RS}).  

In the paper \cite{DK}, Dolgachev and the second author studied the supersingular K3 surface in characteristic 2 with
Artin invariant 1 by using a theory of the Leech lattice.
They showed that the existence of a $(21)_5$-configuration of non-singular  rational curves on the surface and gave a generator of the group of automorphisms of the surface.  Later the authors \cite{Ka} showed the existence of a $(21)_5$-configuration by using a structure of a generalized Kummer surface.  

In this paper, we consider the supersingular K3 surface $X$ with Artin invariant 1
in charactersitic 3.  It is known that $X$ has a rich structure.  For example, $X$ is isomorphic to the Fermat quartic surface on which 112 lines exist, 
the elliptic modular surface with level 4
and the Kummer surface ${\rm Km}(A)$ associated to a superspecial abelian surface $A$ (Shioda \cite{Sh0}, \cite{Sh}).  
Recently Shimada and Zhang \cite{SZ} showed that $X$ is obtained as a purely inseparable triple cover of a non-singular quadric.
Moreover $X$ has extra automorphisms contrary to the case of characteristic zero (see Dolgachev, Keum \cite{DKe}).

The purpose of this paper is to show the existence of 112 non-singular rational curves on $X$ by several ways.
First we shall study $X$ by using a theory of the Leech lattice.
The Picard lattice $S_X$ is isomorphic to the orthogonal complement of the root lattice $A_2\oplus A_2$
in the even unimodular lattice ${\rm II}_{1,25}$ of signature $(1,25)$.  By a description of a fundamental domain of the
reflection subgroup of the orthogonal group ${\rm O}({\rm II}_{1,25})$ due to Conway \cite{Co} and a method of Borcherds \cite{B},
we can see that there exist 112 non-singular rational curves on $X$ (Lemma \ref{112}).  
We show that
these 112 non-singular rational curves are lines on the Fermat quartic surface (\S \ref{fermat}).  
Next we show that $X$ is obtained as a double cover
of a non-singular quadric, and describe 112 curves in terms of non-singular curves on the quadric (\S \ref{double}).  Finally 
we construct 80 elliptic curves and 16 non-singular curves of genus 4 on $A$ whose images, 
together with 16 non-singular curves obtained as exceptional curves, give
112 non-singular rational curves on ${\rm Km}(A)$ (\S \ref{kummer}). 
To construct these 80 elliptic curves, we consider decompositions of the curves of genus 4 into elliptic curves. More precisely, a curve $C$ of genus 4 on $A$ corresponds, 
as a divisor, to a quaternion hermitian $2\times 2$ matrix $M$. We first decompose $M$ into two quaternion hermitian matrices which correspond to elliptic curves $\Delta$ and 
$\Delta'$ respectively.
Namely, we have a decomposition $C \equiv \Delta + \Delta'$ in the N\'eron-Severi
group $NS(A)$. Since we have 10 different decompositions of $M$, we get 20 elliptic curves in this way. By translations by 2-torsion points of $A$, we have 80 elliptic curves
in total. By our construction we can controle the intersection numbers of the elliptic curves and the curves of genus 4, and can show that the elliptic curves have the required
properties.
We will also show that
there exist two sets of non-singular rational curves on $X$ each of which consist of 16 disjoint curves such that
each curve from one set meets exactly 10 members in other set.  Such a pair of two sets is called a
$(16)_{10}$-configuration.  It is known that there exists a $(16)_{10}$-configuration on a Kummer surface in characteristic 0 
(Barth, Nieto \cite{BN}, Naruki \cite{Na}, Traynard \cite{T}).  And we will show that there exist 280 ${\bf F}_9$-rational points on $X$
such that each of 112 curves contains 10 ${\bf F}_9$-rational points and each ${\bf F}_9$-rational point lies on 4 curves.  Such relation is
called a $(280_{4}, 112_{10})$-configuration (for the notation, see Dolgachev \cite{D} for instance).

We thank Daniel Allcock, Igor Dolgachev and JongHae Keum for useful discussions.

\section{Leech roots}\label{leech}

A  lattice $(L, \langle, \rangle)$ is a pair of a free ${\bf Z}$-module $L$ of rank $r$ and
a non-degenerate symmetric integral bilinear form $\langle, \rangle : L \times L \to {\bf Z}$. 
For simplicity we omit $\langle, \rangle$ if there are no confusion.  
We denote by $U$ the hyperbolic plane, that is,
an even unimodular lattice of signature $(1,1)$, and by $A_m, \ D_n$ or $\ E_k$ the even negative definite lattice defined by
the Cartan matrix of type $A_m, \ D_n$ or $\ E_k$ respectively.  We denote by $L\oplus M$ the orthogonal direct sum of lattices $L$ and $M$.
Also we denote by $L^{m}$ the orthogonal direct sum of {\it m}-copies of $L$.
Let $L$ be an even lattice and let $L^* ={\rm Hom}(L,{\bf Z})$.  We denote by $A_L$ the quotient
$L^*/L$ and define a map
$$q_L : A_L \to {\bf Q}/2{\bf Z}$$
by $q_L(x+L) = \langle x, x\rangle\ {\rm mod}\ 2{\bf Z}$.  We call $q_L$ the  discriminant quadratic form.  
Let ${\rm O}(L)$ be the orthogonal group of $L$, that is, the group of isomorphisms of $L$ preserving the bilinear form.
Similarly ${\rm O}(q_L)$ denotes the group of isomorphisms of $A_L$ preserving $q_L$.
Note that there is a natural map from ${\rm O}(L)$ to ${\rm O}(q_L)$.

Let ${\rm II}_{1,25}$ be an even unimodular lattice of signature $(1,25)$.  
A vector $r$ in ${\rm II}_{1,25}$ with $\langle r,r\rangle = -2$ is called a root of ${\rm II}_{1,25}$.  
For each root $r$ we associate the reflection $s_r$ defined by
$$s_r(x) = x +\langle r, x\rangle r$$
which is an isometry of ${\rm II}_{1,25}$.  
Let $W({\rm II}_{1,25})$ be the reflection group of ${\rm II}_{1,25}$, that is, the subgroup of the orthogonal group ${\rm O}({\rm II}_{1,25})$ of 
${\rm II}_{1,25}$  generated by
all reflections.  Obviously $W({\rm II}_{1,25})$ is normal in ${\rm O}({\rm II}_{1,25})$.  The subgroup $W({\rm II}_{1,25})$ is not finite index in 
${\rm O}({\rm II}_{1,25})$, however,
its fundamental domain can be described explicitly by Conway \cite{Co} as follows.  First
we fix an orthogonal decomposition 
$${\rm II}_{1,25} = U \oplus \Lambda$$
where $U$ is the hyperbolic plane and $\Lambda$ is the Leech lattice, i.e., $\Lambda$ is an 
even negative definite unimodular lattice of rank 24 without roots.

We denote a vector in ${\rm II}_{1,25}$ by $(m,n, \lambda)$, where $m, n$ are integers, $\lambda$ is in $\Lambda$, and the norm is given by
$2mn + \langle \lambda, \lambda \rangle$.  We fix a special vector $\rho=(1, 0, 0)$ called a Weyl vector.
A root $r$ in ${\rm II}_{1,25}$ with $\langle r, \rho \rangle =1$ is called a Leech root.
The set $\Delta({\rm II}_{1,25})$ of all Leech roots bijectively corresponds to the set $\Lambda$ by
\begin{equation}\label{corr}
\Lambda \ni \lambda \longleftrightarrow r=(-1-\langle \lambda, \lambda \rangle /2, 1, \lambda) \in \Delta ({\rm II}_{1,25}).
\end{equation}
Note that if $r, r'$ are Leech roots corresponding to $\lambda, \lambda' \in \Lambda$ respectively, then
\begin{equation}\label{corr2}
(r - r')^2 = (\lambda - \lambda')^2.
\end{equation}
Let $P({\rm II}_{1,25})^+$ be a connected component of the set 
$$\{ x \in {\rm II}_{1,25} \otimes {\bf R} : \ \langle x, x\rangle >0\}.$$
Define 
$$C = \{ x \in P({\rm II}_{1,25})^+ : \ \langle x, r\rangle > 0, r \in \Delta({\rm II}_{1,25}) \}.$$
Then $C$ is a fundamental domain of the reflection group $W({\rm II}_{1,25})$ and ${\rm O}({\rm II}_{1,25})$ is a split extension of 
$W({\rm II}_{1,25})$ by ${\rm Aut}(C)$ so that
${\rm Aut}(C)$ is isomorphic to a split extension of 
${\bf Z}^{24}$ by ${\rm O}(\Lambda)$ (Conway \cite{Co}).

Recall that $\Lambda$ is realized as a certain subgroup  in  ${\bf R}^{24} =
{\bf R}^{{\bf P}^1({\bf F}_{23})}$ equipped with inner product $\langle x,y \rangle
= -\frac{x\cdot y}{8}$ (in group theory one often changes the sign to the opposite).  For any subset $A$ of $\Omega = {\bf P}^1({\bf F}_{23})$
let $\nu_A$ denote the vector $\sum_{i\in A} e_i$, where $\{e_\infty,
e_0,\ldots,e_{22}\}$ is the standard basis in ${\bf R}^{24}$. A Steiner
system $S(5,8,24)$ is a set consisting of eight-element subsets of
$\Omega$ such that
any five-element subset of $\Omega$ belongs to a unique element of $S(5,8,24)$. An
eight-element subset in $S(5,8,24)$ is called an octad. Then
$\Lambda$ is
defined as a lattice generated by the vectors $\nu_\Omega-4\nu_\infty$ and
$2\nu_K$, where $K$ belongs to the
Steiner system $S(5,8,24)$.  For more details, see \cite{Co}.

Now we consider the following vectors in $\Lambda$: 
\begin{equation}\label{XYZ}
X = 4\nu_{\infty} + \nu_{\Omega}, Z = 0,
P = 4\nu_{\infty} + 4\nu_{0}, Q = \nu_{\Omega} - 4\nu_{1}.
\end{equation}
Then the corresponding Leech roots
\begin{equation}\label{xyz}
x = (2, 1, X), z = (-1, 1, 0), p = (1, 1, P), q = (1,1,Q)
\end{equation}
genetare a root lattice $R$ isomorphic to $A_2\oplus A_2$.  Obviously $R$ is primitive in ${\rm II}_{1,25}$.   

Let $S$ be the orthogonal complement of $R$ in ${\rm II}_{1,25}$.   
Then the signature of $S$ is $(1, 21)$ and $S^*/S \cong R^*/R \cong ({\bf Z}/3{\bf Z})^2$.
Under the embedding of $S$ in ${\rm II}_{1,25}$ we define
$D=C \cap P(S)^+$ where $P(S)^+$ is a connected component of the set $\{ x \in S \otimes {\bf R} :  \langle x, x\rangle >0\}.$
Then $D$ is non-empty, and has only a finite number of faces and of finite volume (\cite{B}, Lemmas 4.1, 4.2, 4.3).
We denote by ${\rm Aut}(D)$ the group of isometries of $S$ which preserves $D$.  
It is easy to see that the group ${\rm O}(q_R)$ is isomorphic to the dihedral group of order 8 which is represented by
the symmetry group of the Dynkin diagram of $R$.  This implies that  
the natural map
${\rm O}(R) \to {\rm O}(q_R)$ is surjective.  It follows from \cite{N}, Proposition 1.6.1 that 
any isometry of $S$ can be extended to an isometry of ${\rm II}_{1,25}$.
Hence any isometry in ${\rm Aut}(D)$ can be extended to the one in ${\rm O}({\rm II}_{1,25})$ preserving $C$.  
Therefore ${\rm Aut}(D)$ is isomorphic to
the subgroup of ${\rm Aut}(C)$ preserving $R$.

\begin{lm}\label{aut}
{\it The group ${\rm Aut}(D)$ is an extension of ${\rm PSU}(4,{\bf F}_3)$ by a Dihedral group of order $8$ in which
${\rm PSU}(4,{\bf F}_3)$ acts on $R$ trivially and the quotient by it is the symmetry group of the Dynkin diagram of $R$.}
\end{lm}
\proof{See \cite{C}, p.52.}

\begin{lm}\label{112}
{\it There are exactly $112$ Leech roots which are orthogonal to $R$, and ${\rm PSU}(4, {\bf F}_3)$ acts transitively on these Leech roots.}
\end{lm}
\proof{The desired Leech roots correspond to the following norm $(-4)$-vectors in $\Lambda$:

(i)\ 56 $(-4)$-vectors
$2\nu_{K}$ where $K$ is an octad which contains $\infty, 0$, but does not contain $1$;

(ii)\ 56 $(-4)$-vectors
$(x_{\infty}, x_{1},x_{k_{1}},...,x_{k_{6}},...)
= (3,-1,...-1,1,...,1)$

\noindent
where $\{ \infty, 1, k_{1},.., k_{6} \}$ is an octad which does not contain $0$.

On the other hand, the stabilizer of such a vector has index at most 112 in ${\rm Aut}(D)$.  It follows from \cite{C}, p.52,
the table of maximal subgroups of ${\rm PSU}_4(3)$ that the stabilizer is $3^4\cdot {\cal A}_6$ where ${\cal A}_6$ is an alternating group of degree
6.  Thus ${\rm PSU}(4, {\bf F}_3)$ acts transitively on 112 roots.
}

\begin{remark}\label{todd} 
In his paper {\rm \cite{To}}, Todd gave the $759$ octads of  the Steiner system $S(5,8,24)$ 
 {\rm (p.219, Table I).}
The octads in the proof of Lemma {\rm \ref{112}} correspond to the following octads in Todd's table $:$
the last $16$ in the first column, the second and third column.
\end{remark}

\begin{remark}\label{incident} 
The incident relations between the above $112$ Leech roots are as follows
$:$ Recall that there are two types {\rm (i)} and {\rm (ii)} of Leech roots 
$($see the proof of Lemma {\rm \ref{112}}$)$.   
If $r$ and $r'$ are the same types corresponding to octads 
$K$, $K'$, respectively, then $\langle r, r'\rangle = 1$
$($resp. $\langle r, r'\rangle = 0$$)$ iff
$| K \cap K'| = 2$ $($resp. $|K \cap K'| = 4$$)$.  
If $r$ and $r'$ are different types corresponding to octads $K, K'$, respectively, then  
$\langle r, r'\rangle = 1$ $($resp. $\langle r, r'\rangle = 0$$)$ iff
$| K \cap K'| = 4$ $($resp. $|K \cap K'| = 2$$)$.  
\end{remark}

Next we shall classify the hyperplanes bounding $D$.  To do this it is enough to consider the Leech roots $r$ such that
$r$ and $R$ generate a root lattice $R'$ of rank 5 because otherwise the hyperplane $r^{\perp}$ in $P({\rm II}_{1,25})^+$ does not meet
$P(S)^+$.   Obviously there are three cases occur, that is, 
$$R' = A_2\oplus A_2\oplus A_1, \ A_3\oplus A_2, \  A_5.$$
The case $R' = A_2\oplus A_2\oplus A_1$ has been determined in Lemma \ref{112}.

\begin{lm}\label{a5}
{\it There are exactly $5184$ Leech roots $r$ such that $r$ and $R$ generate a root lattice $A_5$, and 
${\rm Aut}(D)$ acts transitively on these Leech roots.
The projection $r'$ of $r$ into $S\otimes {\bf Q}$ has norm $-2/3$.}
\end{lm}
\proof{By considering the action of the dihedral group on the simple roots of $R$, we may assume that $r$ meets $x$ and $p$.
The desired Leech roots correspond to the following norm $(-4)$-vectors in $\Lambda$:

120 $(-4)$-vectors
$2\nu_{K}$ where $K$ is an octad which contains $0$, but does not contain $\infty, 1$;

$120\times 7 = 840$ $(-4)$-vectors
$(x_{k_1}, x_1, x_{k_{2}},...,x_{k_{7}},...) = (3,-1,...,-1,1,...,1)$

\noindent
where $\{k_1, 1, k_{2},.., k_{7} \}$ is an octad which does not contain $\infty, 0$;

$(77-21)\times 6=336$ $(-4)$-vectors 

$(x_1, x_{k_{1}}, x_{\infty}, x_{k_2}, ...,x_{k_{6}},...)
= (-2,-2,2,2,2,2,2,2,0,...,0)$

\noindent
where $\{\infty, 1, k_{1},.., k_{6} \}$ is an octad which does not contain $0$.
Therefore the total number is $1296 \times 4 = 5184$. 
On the other hand, the stabilizer subgroup of $A_5$ 
in ${\rm PSU}(4,{\bf F}_3)$ is a subgroup of $M_{22}$ which fixes an octad.  Note that it fixes one point in the octad and one point 
of the complement of the octad.  Hence the stabilizer is the alternating group ${\cal A}_7$ of degree 7.  The index of ${\cal A}_7$ in
${\rm PSU}(4,{\bf F}_3)$ is $1296$.  Hence the second assertion follows.  
Finally write $r=r'+r'', \ r' \in S\otimes {\bf Q}, r''\in R\otimes {\bf Q}$.
Then we can see $(r'')^2 = -4/3$.  Hence $(r')^2=-2/3$.
}

\begin{lm}\label{a3}
{\it There are exactly $648$ Leech roots $r$ such that $r$ and $R$ generate $A_3\oplus A_2$, and 
${\rm Aut}(D)$ acts transitively on these Leech roots.  The projection $r'$ of $r$ into $S\otimes {\bf Q}$ has norm $-4/3$.}
\end{lm}
\proof{By considering the action of the dihedral group on the simple roots of $R$,
we may assume that $r$ meets $q$.
Then the desired Leech roots correspond to the following norm $(-4)$-vectors in $\Lambda$:

120 $(-4)$-vectors
$2\nu_{K}$ 

\noindent
where $K$ is an octad which contains $\infty$, but does not contain $0, 1$;

21 $(-4)$-vectors
$(x_{\infty}, x_0, x_{1},x_{k_{1}},...,x_{k_{5}},...)
= (3,-1,...-1,1,...,1)$

\noindent
where $\{ \infty, 0, 1, k_{1},.., k_{5} \}$ is an octad;

21 $(-4)$-vectors 
$\nu_{\Omega} - 4\nu_{k}$ 

\noindent
where $k \not= \infty, 0, 1$.  
Therefore the total number is $162 \times 4 = 648$.  On the other hand, the stabilizer subgroup of $A_3\oplus A_2$ 
in ${\rm PSU}(4,{\bf F}_3)$ is a subgroup of $M_{24}$ which fixes three points of 24 letters, that is, $L_3(4)$.  The index of $L_3(4)$ in
${\rm PSU}(4,{\bf F}_3)$ is $162$.  Hence the second assertion follows.  
Finally write $r=r'+r'', \ r' \in S\otimes {\bf Q}, r''\in R\otimes {\bf Q}$.
Then we can see $(r'')^2 = -2/3$.  Hence $(r')^2=-4/3$.
}


\begin{lm}\label{weyl}
{\it Let $w$ be the projection of the Weyl vector $\rho = (1,0,0) \in U \oplus \Lambda$ into $S^*$.  Then
$w \in S$ and $\langle w, w\rangle = 4$.  Moreover $w$ is contained in $D$.  In particular, there are no $(-2)$-vectors in $S$
perpendicular to $w$.}
\end{lm}
\proof{Let $\rho = w + w'$ where $w' \in R^*$.  Recall that $x, z, p, q$ generate $R$.  Since $\langle x, w\rangle 
=\langle z, w'\rangle = \langle p, w'\rangle = \langle q, w'\rangle =1$,
we have $w' = -(x+z+p+q)$.  Hence $w' \in R$.
This implies that $w \in S$ and $\langle w,w\rangle = -\langle w',w'\rangle = 4$.  
It is known that $w$ is contained in $D$ (\cite{B}).}

\begin{lm}\label{picard}
Denote by $\{ r_{\alpha}\}$ the set of all $112$ Leech roots in Lemma {\rm \ref{112}}.  Then

$(1)$ Each $r_{\alpha}$ meets exactly $30$ members from $\{r_{\alpha}\}$.

$(2)$ $\{ r_{\alpha}\}$ generate $S$. 

$(3)$   $\langle w, r_{\alpha}\rangle =1$.

$(4)$ $w = {1\over 28} \sum_{\alpha} r_{\alpha}.$
\end{lm}
\proof{
(1) Since ${\rm PSU}(4,{\bf F}_3)$ acts on the set of 112 Leech roots transitively, it is enough to see the assertion for
one of them.  We take the Leech root $r_K$ corresponding to $2\nu_K$ where
$$K=\{\infty, 0, 2,3,4,8,9,21\}$$
(see the proof of Lemma \ref{112}).   By using the incident relations mentioned in Remark \ref{incident}, we can easily see that the number of Leech roots $r_{\alpha}$ with $\langle r_K, r_{\alpha}\rangle =1$
(resp. $\langle r_K, r_{\alpha}\rangle = 0$) is 30 (resp. 81).  The Leech roots $r_{\alpha}$ with $\langle r_K, r_{\alpha} \rangle = 1$ correspond to
the following 30 octads :
$$K_{1} = \{ \infty, 0, 5,6,7,13,16,17 \}, \quad K_{2} = \{ \infty, 0,5,7,11,14,18,19 \}, \\$$
$$K_{3} = \{ \infty, 0,5,10,13,14,15,22 \}, \quad K_{4} = \{ \infty, 0,5,11,12,15,17,20 \},\\$$
$$K_{5} =\{ \infty, 0, 6,7,10,12,15,18 \}, \quad K_{6} = \{ \infty, 0,6,10,14,17,19,20 \},\\$$
$$ K_{7} = \{ \infty, 0,6,11,12,13,19,22 \},  \quad K_{8} = \{ \infty, 0,7,15,16,19,20,22 \}, \\$$
$$K_{9} = \{ \infty, 0, 10,11,13,16,18,20 \}, \quad K_{10} = \{ \infty, 0, 12,14,16,17,18,22 \},\\$$
$$K'_{1} = \{ \infty, 1,2,5,6,8,9,16 \}, \quad K'_{2} = \{ \infty, 1,2,4,7,9,11,14 \},\\$$
$$K'_{3} = \{ \infty, 1,2,8,13,14,15,21 \}, \quad K'_{4} = \{ \infty, 1,2,4,5,12,20,21 \},\\$$
$$K'_{5} = \{  \infty, 1,2,3,4,6,15,18 \}, \quad K'_{6} = \{ \infty, 1,2,4,8,10,17,19 \}, \\$$
$$ K'_{7} = \{ \infty, 1,2,3,9,12,13,19 \}, \quad K'_{8} = \{ \infty, 1,2,3,7,8,20,22 \},\\$$
$$K'_{9} = \{ \infty, 1,2,3,10,11,16,21 \}, \quad K'_{10} = \{ \infty,1,2,9,17,18,21,22 \}, \\ $$
$$K''_{1} = \{ \infty, 1,3,4,5,9,10,22 \}, \quad K''_{2} = \{ \infty, 1,3,5,8,18,19,21 \},\\$$
$$K''_{3} = \{ \infty, 1,3,4,5,9,10,22 \}, \quad K''_{4} = \{ \infty, 1,3,8,9,11,15,17 \},\\$$
$$K''_{5} = \{ \infty, 1,7,8,9,10,12,21\}, \quad K''_{6} = \{ \infty, 1,3,6,9,14,20,21 \},\\$$
$$K''_{7} = \{ \infty, 1,4,6,8,11,21,22 \}, \quad K''_{8} = \{ \infty, 1,4,9,15,16,19,21 \}, $$
$$K''_{9} = \{ \infty, 1,4,8,9,13,18,20 \},\quad K''_{10} = \{ \infty, 1,3,4,8,12,14,16 \}.\\$$

(2) Denote by $k_i, k'_j, k''_l$ the corresponding Leech roots to $K_i, K'_j, K''_l$ respectively.  
Then we can easily see that 
$\langle k_i, k'_i\rangle = \langle k'_i, k''_i \rangle = \langle k''_i, k_i\rangle = 1$ $(1\leq i \leq 10)$ 
and all other intersection numbers between $k_i, k'_j, k''_l$ are zero.  
Therefore $(k_i + k'_i + k''_i)^2 =0$ and $\langle k_i+k'_i+k''_i, k_j+k'_j+k''_j\rangle = 0$ for $i\not= j$.
Hence $ k_i+k'_i+k''_i = k_j+k'_j+k''_j$.  The remaining $81 (= 112 - 31)$ Leech roots have non-zero intersection number with 
$ k_i+k'_i+k''_i$.  Since each of $k_i, k'_j, k''_l$ has intersection number 1 with $27 (=30-3)$ members among the remaining $81$ Leech roots.
Hence each of the remaining 81 Leech roots has intersection number 1 with $10 = 30 \times 27 / 81$ members among $k_i, k'_j, k''_l$.
Thus if $f$ is any one of the remaining $81$ Leech roots, then $\langle f,  k_i+k'_i+k''_i\rangle = 1$.
Then two vectors $f$ and $k_i+k'_i + k''_i$ generate a hyperbolic plane $U$ in $S$, and
its orthogonal complement contains $A_2^{10}$ generated by $k_i, k'_i, k''_i$ perpendicular to $f$.
Let $S'$ be the sublattice of $S$ generated
by $31$ Leech roots $k_i, k'_j, k''_l$ and $f$.  Then $S' \cong U\oplus A_2^{10}$, $S'$ is of finite index in $S$ and the quotient $S/S'$ is 
isomorphic to $({\bf Z}/3{\bf Z})^4$.  The other Leech roots except $r_K, k_i, k'_j, k''_l$ bijectively correspond to $({\bf Z}/3{\bf Z})^4$.
Therefore $112$ Leech roots generate a sublattice of $S$ with the discriminant
$\pm 3^{10}/{81}^2 = \pm 3^2$.  Since $S$ has the discriminant $\pm 3^2$, 112 Leech roots generate $S$.  
In the following Lemma \ref{a2fibration}, we shall give a geometric interpretation of this argument.

(3) Since $\langle \rho, r_{\alpha} \rangle = 1$, we have $\langle w, r_{\alpha} \rangle = 1$.  

(4)  The right hand side $\sum r_{\alpha} /28$ 
has the same property in (3).   Since $112$ Leech roots generate $S$, the assertion follows.}

The following shows the existence of a $(16)_{10}$-configuration among 112 Leech roots.  Later we shall show the existence 
of $(16)_{10}$-configuration geometrically (see Theorem \ref{16-10th}).

\begin{prop}\label{16-10}
{\it Among $112$ Leech roots, there exist two sets ${\cal A}$ and ${\cal B}$ each of which consists of $16$ Leech roots such that each one from
one set meets exactly $10$ members of the other set.}
\end{prop}
\proof{
We use the same notation as in the proof of Lemma \ref{112}.  Let ${\cal A}$ be the set of Leech roots corresponding
$2\nu_K$ where $K$ is an octad containing $\infty, 0, 2$, but does not contain $1$.  And let ${\cal B}$ be the set of Leech roots
correspond to $(-4)$-vectors
$(x_{\infty}, x_{1},x_{k_{1}},...,x_{k_{6}},...)
= (3,-1,...-1,1,...,1)$
where $\{ \infty, 1, k_{1},.., k_{6} \}$ is an octad which contains 2 and does not contain $0$.  We can easily see that  
the pair of ${\cal A}$ and $\cal B$ is a $(16)_{10}$-configuration.
 }

\begin{lm}\label{picard2}
{\it Let $X$ be the supersingular $K3$ surface in characteristic $3$ with Artin invariant $1$.  
Let $S_X$ be the Picard lattice of $X$.  Then $S_X$ is isomorphic to $S$.}
\end{lm}
\proof{
It is known that $S_X \cong U\oplus E_8\oplus E_8\oplus A_2\oplus A_2$ (see \cite{RS},\ \S 5).
Since $S$ is the orthogonal complement of $R = A_2\oplus A_2$ in the unimodular lattice ${\rm II}_{1,25}$,
$q_R = -q_S$ (Corollary 1.6.3 in \cite{N}).  By using these facts, we have $q_{S_X} \cong q_R \cong -q_R \cong q_S$.
It now follows from Theorem 1.14.2 in \cite{N} that
$S$ is isomorphic to $S_X$.   
 }

Let $P(S_{X})^+$ be the connected component of the cone $\{ x \in S_{X} \otimes
{\bf R}  :  \langle x, x \rangle > 0 \}$ which contains an ample divisor.
Let $D(X)$ be the ample cone of $X$ which consist of all vectors $x$ in $P(S_X)^+$ satisfying
$\langle x, \delta \rangle > 0$ for any effective divisor $\delta$ on $X$ with $\langle \delta, \delta \rangle = -2$.
Under an identification of $S_X$ with $S$, we may assume that $P(S_X)^+$ is contained in $P({\rm II}_{1,25})^+$.
Since any $(-2)$-vector in $S_X$ is a $(-2)$-vector in ${\rm II}_{1,25}$, we may assume that $D = C \cap P(S_X)^+$ is contained in $D(X)$.
By Lemma \ref{weyl}, the projection $w$ of the Weyl vector $\rho$ is contained in $D(X)$, that is, $w$ is an ample class.

\begin{remark}\label{auto} 
In case of the supersingular $K3$ surface $Y$ in characteristic $2$ with Artin invariant $1$, the group ${\rm Aut}(Y)$ of automorphisms of $Y$  
is generated by
${\rm PGL}(3, {\bf F}_4)\cdot {\bf Z}/2{\bf Z}$ and $168$ involutions $($see {\rm \cite{DK}}$)$.  
The restriction of the Conway's fundamental domain $C$ 
to the positive cone $P(S_Y)^+$ of $Y$ consists of $42$ faces defined by $42$ $(-2)$-vectors and $168$ faces defined by $168$ $(-1)$-vectors.
The group ${\rm PGL}(3, {\bf F}_4)\cdot {\bf Z}/2{\bf Z}$
preserves the $42$ $(-2)$-vectors, and $168$ involution act on $S_Y$ as reflections with respect to $168$ $(-1)$-vectors.
It would be interesting to give a generator of ${\rm Aut}(X)$ in our case.
\end{remark}

\section{The Fermat quartic surface}\label{fermat}

It is known that the Fermat quartic surface in ${\bf P}^3$ defined by
\begin{equation}\label{fermat1}
x_0^4+x_1^4+x_2^4+x_3^4 = 0
\end{equation}
is a supersingular
$K3$ surface with Artin invariant 1 in characteristic 3 (\cite{Sh}), and hence is isomorphic to $X$.
In the following we identify $X$ with the Fermat quartic surface.
Note that the equation (\ref{fermat1}) is a hermitian form over ${\bf F}_9$, and hence the group
${\rm PGU}(4,{\bf F}_3)$ acts on $X$.   Since ${\rm PSU}(4, {\bf F}_3)$ is simple, it acts 
trivially on $S_X^*/S_X$.  Therefore under the identification of $S_X$ with $S$
(Lemma \ref{picard2}), the action of
${\rm PSU}(4, {\bf F}_3)$ on $S$ can be extended to the one on $U\oplus \Lambda$ acting
trivially on $R=A_2\oplus A_2$.  It follows from Lemma \ref{aut} that ${\rm PSU}(4, {\bf F}_3)$ is isomorphic to the point wise stabilizer group of $R$.  It is known that the smallest dimension of irreducible representations 
of ${\rm PSU}(4, {\bf F}_3)$ is 21 (see \cite{C}, page 54).  Thus the projection $w$ of the Weyl vector 
$\rho$ is a unique invariant vector with norm $4$ in $S_X$ under the action of ${\rm PSU}(4, {\bf F}_3)$, and hence $w$
coincides with the class of the hyperplane section of the Fermat quartic surface.
By Lemma \ref{picard}, (3), we conclude that 112 Leech roots in Lemma \ref{112} are lines on $X$.
We remark that it is known that the Fermat quartic surface contains exactly 112 lines (\cite{Se2}).

Recall that the number of isotropic vectors of a $n$-dimensional non-degenerate hermitian form over ${\bf F}_{q^2}$ is given by
$$(q^n - (-1)^n)(q^{n-1} - (-1)^{n-1}) + 1.$$
Therefore the number of ${\bf F}_{9}$-rational points on the Fermat quartic surface $X$ is ${(3^4 - 1)(3^3 +1) \over (3^2 -1)} = 280$.  
On the other hand, ${\bf F}_9$-rational points of ${\bf P}^1$ is $10$.  Therefore ${112 \times 10 \over 280} = 4$ lines pass at each 
${\bf F}_9$-rational point of $X$.  This relation between 112 lines and 280 rational points is called $(280_4, 112_{10})$-configuration.
For the notion of configuration, we refer the reader to \cite{D}.

\begin{prop}\label{280-4}
{\it There exists a $(280_4, 112_{10})$-configuration on $X$}
\end{prop}
Later we shall give this configuration in terms of the Kummer surface (see Theorem \ref{280-4th}).
The following Lemma is a geometric interpretation of Lemma \ref{picard}.

\begin{th}\label{a2fibration}
{\it For each line $l$ on $X$, there exists exactly $30$ lines meeting with $l$.  These $30$ lines are components of $10$ singular fibers of
type {\rm IV} of
a quas-elliptic fibration, and $l$ is a cuspidal curve of this fibration.
The other $81$ lines are sections of this fibration.  The linear system of hyperplane sections of the Fermat quartic surface containing $l$ defines this fibration.
}\end{th}
\proof{
We use the same notation as in the proof of Lemma \ref{picard}.
The linear system $| k_i + k'_i + k''_i|$ defines a quas-elliptic or an elliptic fibration $\pi : X \to {\bf P}^1$.
Obviously $\pi$ has 10 singular fibers of type $\tilde{A}_2$ and 81 sections.  It now follows from \cite{RS}, Theorem in \S 4
that $\pi$ is a quasi-elliptic fibration.   Recall that a singular fiber of type ${\rm IV}$ of a quasi-elliptic fibration consist of three components meeting at one point and the cuspidal curve passes through singular points on each fiber.  Since the cuspidal curve and the line $l$ have the same intersection numbers with components of fibers and sections.  Hence the line $l$ coincides the cuspidal curve.
}

\section{A double quadric model}\label{double}

In this section, we shall construct $X$ as a double cover of a non-singular quadric surface. 
We denote by $\zeta$ a primitive eighth root of unity satisfying $\zeta^2 + \zeta =1$ and $\zeta^2 =\sqrt{-1}$.
Let $Q$ be a non-singular quadric ${\bf P}^1 \times {\bf P}^1$ and let $((u_0:u_1),(v_0:v_1))$ be a homogenous bi-coordinate of $Q$.
Consider non-singular rational curves $C, C', C_i, D_i$ defined by
$$C: \ u_0v_0^3 = u_1v_1^3, \ \ C' : \ u_0^3v_0=u_1^3v_1;$$
$$C_{i}: u_1=\zeta^iu_0 \ (1\leq i \leq 8), \ C_9: u_1=0, \ C_{10} : u_0=0;$$
$$D_{i}: v_1=\zeta^iv_0 \ (1\leq i \leq 8), \ D_9: v_1=0, \ D_{10} : v_0=0.$$
Then $C$ and $C'$ meet at 10 points $((0:1),(1:0)), ((1:0),(0:1)), ((1:\zeta^i), (1:\zeta^{5i})) \ (1\leq i\leq 8)$
transversally.  Note that $C_i$ (resp. $D_j$) meets $C$ (resp. $C'$) at one of the above 10 points with multiplicity 3.
Also there are 30 curves of bidegree $(1,1)$ passing through 4 points from the above 10 points.
They are defined by:
$$ u_0v_0 - \zeta^{2k} u_1v_1=0 \ (0 \leq k \leq 3); \ u_0v_1 \pm u_1v_0=0;$$
$$u_0v_0 \pm u_0v_1+u_1v_1 =0; \ u_0v_0 \pm \zeta u_0v_1- \zeta^2u_1v_1 =0;$$
$$ \ u_0v_0 \pm \zeta^3 u_0v_1+\zeta^2u_1v_1 =0; 
\ u_0v_0 \pm \zeta^2 u_0v_1-u_1v_1 =0;$$
$$u_0v_0 \pm u_1v_0+u_1v_1 =0; \ u_0v_0 \pm \zeta u_1v_0- \zeta^2u_1v_1 =0;$$
$$ \ u_0v_0 \pm \zeta^3 u_1v_0+\zeta^2u_1v_1 =0; \ u_0v_0 \pm \zeta^2 u_1v_0-u_1v_1 =0;$$
$$u_0v_0 + \zeta^{2k}(u_0v_1 -u_1v_0) -u_1v_1=0 \ (0 \leq k \leq 3);$$
$$u_0v_0 \pm \zeta (u_0v_1+u_1v_0) -\zeta^2 u_1v_1 =0;  \ u_0v_0 \pm \zeta^3(u_0v_1+u_1v_0) + \zeta^2u_1v_1 =0.$$
Let $Z$ be the double cover of $Q$ branched along the divisor $C+C'$ of bidegree $(4,4)$.  Then $Z$ has 10
nodes over the above 10 points.  Let $\tilde{Z}$ be the minimal resolution.  Denote by $E_1,..., E_{10}$ the
exceptional curves, and $l, l'$ the inverse images of $C, C'$ respectively.  
Let $C_i^{\pm}$ (resp. $D_j^{\pm}$)
be the inverse images of $C_i$ (resp. $D_j$).  Also the 30 curves of bidegree $(1,1)$ as above split 60 non-singular rational curves on $\tilde{Z}$.
Then $| C_i^+ + C_i^- + E_i|$ defines a quasi-elliptic fibration $\pi$ with ten singular fibers
of type IV. 
This fibration $\pi$ is nothing but the one induced from the projection from $Q$ to the first factor ${\bf P}^1$.  
The cuspidal curve $l$ passes through the singular points of each singular fiber.  The fibration $\pi$ has 81 sections
$l'$, $D_j^{\pm}$ ($1\leq j \leq 10$), and the above 60 non-singular rational curves.  
We now conclude that the rank of Picard lattice is 
$22$ and its discriminant is $\pm 3^{10}/(81)^2 = \pm 3^2$.  Thus we have the following theorem.

\begin{th}\label{}
{\it  $\tilde{Z}$ is the supersingular $K3$ surface with the Artin invariant $1$.
In particular $\tilde{Z}$ is isomorphic to $X$.
}\end{th}

\begin{remark}\label{segre}
The covering transformation of $\tilde{Z} \to Q$ is the Segre involution with respect to a pair of two skew lines $l$, $l'$ 
$($see \rm{\cite{Se1}}, \S $16)$.
\end{remark}

\section{The N\'eron-Severi group of a superspecial abelian surface}
Let $k$ be an algebraically closed field of characteristic $p > 0$
and let $E$ be a supersingular elliptic curve defined over $k$.
We consider a superspecial abelian surface $A = E_{1} \times E_{2}$
with $E_{1} = E_{2} = E$. We denote by $E_{1}$ (resp. $E_{2}$) 
the divisor $E_{1}\times \{0\}$ (resp. the divisor $\{0\}\times E_{2}$)
on $A$ if we have no confusion. We set $X = E_{1} + E_{2}$. 
$X$ is a principal polarization on $A$.
We set ${\cal O} = {\rm End}(E)$ and $B = {\rm End}^{0}(E)=
{\rm End}(E)\otimes {\bf Q}$.
Then, $B$ is a quaternion division algebra over the rational number field ${\bf Q}$
with discriminant $p$, and ${\cal O}$ is a maximal order of $B$.
For a divisor $L$, we have a homomorphism
$$
\begin{array}{cccc}
   \varphi_{L} :  &A  &\longrightarrow  &{\rm Pic}^{0}(A) \\
      & x & \mapsto & T_{x}^{*}L - L,
\end{array}
$$
where $T_{x}$ is the translation by $x \in A$.
We set 
$$
     H = \{
\left(
\begin{array}{cc}
\alpha  & \beta \\
\gamma & \delta
\end{array}
\right)
~\mid
~\alpha, \delta \in {\bf Z},~\gamma = \bar{\beta}
\}.
$$
Then we have the following theorem (cf. \cite{K1} and \cite{K2}).
\begin{th}\label{intersection}
The homomorphism
$$
\begin{array}{cccc}
j : & NS(A) & \longrightarrow  & H \\
   & L  &\mapsto  & \varphi_{X}^{-1}\circ\varphi_{L}
\end{array}
$$
is bijective.  By this correspondence, we have
$$
j(E_{1}) = \left(
\begin{array}{cc}
0  & 0 \\
0 & 1
\end{array}
\right),~
j(E_{2}) = \left(
\begin{array}{cc}
1  & 0 \\
0 &  0
\end{array}
\right).
$$
For $L_{1}, L_{2} \in NS(A)$ such that
$$
  j(L_{1}) =
\left(
\begin{array}{cc}
\alpha_{1}  & \beta_{1} \\
\gamma_{1} & \delta_{1}
\end{array}
\right),~
  j(L_{2}) =
\left(
\begin{array}{cc}
\alpha_{2}  & \beta_{2} \\
\gamma_{2} & \delta_{2}
\end{array}
\right),
$$
the intersection number  $(L_{1}, L_{2})$ is given by
$$
(L_{1}, L_{2}) = \alpha_{2}\delta_{1} + \alpha_{1}\delta_{2} - \gamma_{1}\beta_{2}
-\gamma_{2}\beta_{1}.
$$
In particular, for $L\in NS(A)$ such that 
$j(L) = \left(
\begin{array}{cc}
\alpha  & \beta \\
\gamma & \delta
\end{array}
\right)
$
we have
$$
\begin{array}{l}
L^{2} = 2\det \left(
\begin{array}{cc}
\alpha  & \beta \\
\gamma & \delta
\end{array}
\right)\\
(L, E_{1} ) = \alpha,~(L, E_{2}) = \delta .
\end{array}
$$
\end{th}
Let  $m : E\times E \rightarrow E$ be the addition of $E$, and we set
$$
\Delta = {\rm Ker}~ m.
$$
We have $\Delta = \{(P, - P)  ~\mid~P\in E \}$. For two endomorphisms
$a_{1}, a_{2}\in {\rm End}(E)$, we set
$$
\Delta_{a_{1}, a_{2}} = (a_{1}\times a_{2})^{*}\Delta.
$$
Using this notation, we have $\Delta = \Delta_{1,1}$.
We have the following theorem (cf. \cite{K2})
\begin{th}\label{div}
$$
j(\Delta_{a_{1}, a_{2}}) =
\left(
\begin{array}{cc}
\bar{a}_{1}a_{1}  & \bar{a}_{1}a_{2}\\
\bar{a}_{2}a_{1} & \bar{a}_{2}a_{2}
\end{array}
\right).
$$
In particular, we have
$$
j(\Delta) =
\left(
\begin{array}{cc}
1  & 1\\
1 & 1
\end{array}
\right)
$$
\end{th}

\section{Supersingular elliptic curve in characteristic 3}
We summarize, in this section, known facts 
on the supersingular elliptic curve in characteristic 3
which we will use later.  

We have, up to isomorphism, only one supersingular elliptic curve 
defined over $k$, which is given by the equation
$$
               y^{2}  = x^{3} -x .
$$
We denote by $E$ a nonsingular complete model of the supersingular 
elliptic curve, which is defined by
$$
  Y^{2}Z  = X^{3}- XZ^{2},
$$
in the projective plane ${\bf P}^{2}$,
where $(X, Y, Z)$ is a homogeneous coordinate of ${\bf P}^{2}$.
An invariant non-zero regular vector field on $E$ is given by $\frac{\partial}{\partial y}$.
In the affine model, let $(x_{1}, y_{1})$ and $(x_{2}, y_{2})$ be
two points on E. Then, the addition of E is given by
$$
\begin{array}{l}
x = -x_{1} - x_{2} + ((y_{2} - y_{1})/(x_{2} - x_{1}))^{2},\\
y = y_{1} + y_{2} - ((y_{2} - y_{1})/(x_{2} - x_{1}))^{3}.
\end{array}
$$
We denote by $[n]_{E}$ the multiplication by an integer $n$.  $[2]_{E}$ is concretely
given by
$$
\begin{array}{l}
    x = x_{1} + 1/y_{1}^{2},\\
    y = -y_{1} - 1/y_{1}^{3}.
\end{array}
$$
In case of characteristic $p >0$,
we denote by F (resp. V) the relative Frobenius morphism (resp.
the Vershiebung), which has the following relations:
$$
FV = p,~V =\bar{F} = -F,~F^{2} = -p.
$$
We have the following lemma.
\begin{lm}
Let $p$ be a prime number, and let ${\bf F}_{p^{2}}$ be
the finite field with $p^{2}$ elements. For a supersingular elliptic curve $E$
defined over ${\bf F}_{p}$,
we have ${\rm Ker}[p + 1]_{E} = E({\bf F}_{p^{2}})$. In particular, 
$\mid E({\bf F}_{p^{2}})\mid = (p+ 1)^{2}$.
\end{lm}
\proof{A point $P \in E$ is contained in $E({\bf F}_{p^{2}})$ if and only if 
$F^{2}(P) = P$. Since $F^{2} = -p$, we have $F^{2}(P) = P$ if and only if 
$[p + 1]_{E}(P) = 0$.}
Then, we have
$$
\begin{array}{l}
E({\bf F}_{3}) = \{(0, 1,0),(0,0,1), (1, 0,1), (-1, 0, 1)\}, \\
E({\bf F}_{9}) = \{(0,1,0),(0,0,1), (1, 0,1), (-1, 0, 1), (\zeta, \pm \zeta^{3} ,1),
 (\zeta^{2}, \pm \zeta, 1), \\
\hspace*{2.0cm}(\zeta^{3}, \pm \zeta, 1), (\zeta^{5}, \pm \zeta, 1), 
(\zeta^{6}, \pm \zeta^{3},1), (\zeta^{7}, \pm \zeta^{3} ,1)\}.
\end{array}
$$
Here, $\zeta$ is a primitive eighth root of unity which satisfies $\zeta^{2} + \zeta = 1$
and $\zeta^{2} = \sqrt{-1}$.
We set 
$$
P_{\infty} = (0,1,0), P_{0} = (0,0,1), P_{1} = (1,0,1), P_{-1} = (-1, 0, 1).
$$
The point $P_{\infty}$ is the zero point of $E$, and the group of 
4-torsion points of $E$ is given by
$$
 E({\bf F}_{9}).
$$ 
Then $E$ has the following automorphisms $\sigma$ and $\tau$:
$$
\begin{array}{l}
  \sigma : x\mapsto x + 1,~y \mapsto y\\
 \tau : x\mapsto -x ,~y \mapsto \sqrt{-1}y,
\end{array}
$$
which satisfiy
$$
     \sigma^{3} = {\rm id},~ \tau^{2} = -{\rm id},~ 
\tau\circ \sigma = \sigma^{2}\circ \tau.
$$
By this relation, we know that in the maximal order ${\cal O} ={\rm End}(E)$
$$
\bar{\sigma} = \sigma^{2},~ \tau^{2} = -1, ~ \bar{\tau} = -\tau, 
F\circ \sigma = \sigma \circ F, F\circ \tau = - \tau \circ F.
$$
We set $\rho = \tau \circ \sigma$. Then, $\rho$ is of order 4 and we have
$$
\sigma_{*}(\frac{\partial}{\partial y}) = \frac{\partial}{\partial y},~
\tau_{*}(\frac{\partial}{\partial y}) = \sqrt{-1}\frac{\partial}{\partial y},~
\rho_{*}(\frac{\partial}{\partial y}) = \sqrt{-1}\frac{\partial}{\partial y}.
$$
The translation $T_{P_{0}}$ by the point $P_{0}$ is given by
$$
   T^{*}_{P_{0}}x_{1} = -1/x_{1}, ~T^{*}_{P_{0}}y_{1} = y_{1}/x_{1}^{2}.
$$
We set $\theta = T^{*}_{P_{0}}$. Then the group $G$ which is generated by 
$\theta$ acts on the function field $k(E)$, and 
$$
  \zeta^{2}(1/x_{1} -x_{1}), ~ -\zeta (y_{1}/x_{1}^{2} + y_{1}) .
$$
are invariant under the action of $G$. It is easy to see that the invariant field
$k(E)^{G}$ is given by
$$
k(E)^{G} = k(\zeta^{2}(1/x_{1} -x_{1}),  -\zeta (y_{1}/x_{1}^{2} + y_{1})).
$$
We set 
$$
x = \zeta^{2}(1/x_{1} -x_{1}), y = -\zeta (y_{1}/x_{1}^{2} + y_{1}).
$$
Then, $x$ and $y$ satisfy the equation
$$
   y^{2} = x^{3} - x.
$$
Hence, we have $ E/G \cong E$ and we see that the projection
$$
   \pi : E \longrightarrow E
$$
is given by
$$
      x = \zeta^{2}(1/x_{1} -x_{1}), ~y = -\zeta (y_{1}/x_{1}^{2} + y_{1}).
$$
We have
$$
\pi_{*}(\frac{\partial}{\partial y_{1}}) = \zeta \frac{\partial}{\partial y}.
$$
By direct calculations, we have relations:
$$
\begin{array}{l}
\pi = {\rm id} -\tau\\
\tau\circ \pi = \pi \circ\tau = {\rm id} + \tau\\
\sigma\circ \pi = \pi \circ \sigma \\
\pi\circ \pi \circ \tau = \tau \circ\pi\circ \pi  = [2]_{E}\\
{\rm id} + F = -[2]_{E}\circ \sigma
\end{array}
$$
By the general theory of quaternion algebra over ${\bf Q}$, the quaternion algebra
with descriminant 3 is given by
$$
\begin{array}{c}
    B = {\bf Q} + {\bf Q}\alpha + {\bf Q}\beta+ {\bf Q}\alpha\beta \\
    \alpha^{2} = -3,~\beta^{2} = -1,~ \alpha\beta = -\beta\alpha
\end{array}
$$
and  a maximal order ${\cal O}$ of $B$ is given by
$$
{\cal O} = {\bf Z} + {\bf Z}\beta + {\bf Z}(1 + \alpha)/2+ {\bf Z}\beta(1 + \alpha)/2.
$$
(cf. \cite{I} and \cite{K1}). In the case of the supersingular elliptic curve in characteristic 3,
we can take $\alpha$ and $\beta$ as
$$
\beta = \tau, \alpha = F
$$
and then
$$
(1 + \alpha)/2 = -\sigma,~ \beta(1 + \alpha)/2 =-\tau\sigma
$$
by the relations of endomorphisms above.

We consider now in characteristic 3 the superspecial abelian surface 
$A = E_{1} \times E_{2}$ with $E_{1} = E_{2} = E$. Using the notation
in Section 5, on $A$ we take six divisors
$$
E_{1}, ~E_{2},~ \Delta =\Delta_{1, 1}, ~\Delta_{1, \tau},~\Delta_{1, -\sigma},
~\Delta_{1, -\tau\sigma}.
$$
By Theorems~\ref{div}, we have the following lemma.
\begin{lm}\label{basis}
We have the following expressions of divisors:
$$
\begin{array}{l}
j(E_{1})=\left(
\begin{array}{cc}
0  & 0 \\
0 & 1
\end{array}
\right), ~j(E_{2})=\left(
\begin{array}{cc}
1  & 0 \\
0 & 0
\end{array}
\right),~ j(\Delta) = \left(
\begin{array}{cc}
1  & 1 \\
1 & 1
\end{array}
\right),\\
j(\Delta_{1, \tau})=\left(
\begin{array}{cc}
1  & \tau \\
-\tau & 1
\end{array}
\right),
~j(\Delta_{1, -\sigma})=\left(
\begin{array}{cc}
1  & -\sigma \\
-\sigma^{2} & 1
\end{array}
\right),\\
~j(\Delta_{1, -\tau\sigma})=
\left(
\begin{array}{cc}
1  & -\tau\sigma \\
\sigma^{2}\tau & 1
\end{array}
\right).
\end{array}
$$
Moreover, these six divisors give a basis of $NS(A)$.
\end{lm}

By direct calculation based on Theorem~\ref{intersection} we have the following table for the intersection numbers.

\begin{center}
\begin{tabular}{| l |c|c|c|c|c|c|}
\hline
     & $E_{1}$&$ E_{2}$& $\Delta$ &$\Delta_{1, \tau}$&$\Delta_{1, -\sigma}$
& $\Delta_{1, -\tau\sigma} $\\
\hline
$E_{1}$ & 0 & 1 & 1& 1& 1& 1 \\
\hline
$ E_{2}$ & 1 & 0 & 1 & 1 & 1 & 1  \\
\hline
$ \Delta$ & 1 & 1& 0 & 2& 1& 2\\
\hline
$\Delta_{1, \tau}$& 1 & 1& 2 & 0 & 2 & 1\\
\hline
$\Delta_{1, -\sigma}$&  1 & 1& 1& 2& 0 & 2\\
\hline
$\Delta_{1, -\tau\sigma}$& 1 & 1 & 2& 1& 2 & 0\\

\hline
\end{tabular}
\end{center}

\section{A curve of genus 4}
In this section we consider the non-singular complete curve $C$ 
of genus 4 defined by
the equation
$$
    Y^{2} = X^{9} - X.
$$ 
The ${\bf F}_{9}$-rational points are given by
$$
      C( {\bf F}_{9}) = \{\infty, (0,0), (\zeta^{i}, 0)~\mid~i = 0, 1, \ldots, 7\}.
$$
Therefore, $C$ has 10 ${\bf F}_{9}$-rational points.
We have a morphism $\varphi$ from $C$ to $E$ defined by
$$
   x = X^{3} + X, ~ y = Y.
$$
We have
$$
   \varphi_{*}(\frac{\partial}{\partial Y}) = \frac{\partial}{\partial y}
$$ 
We need to construct one more morphism from $C$ to $E$. For this purpose,
we consider an automorphism $\eta$ of $C$ which is defined by
$$
       X \mapsto (X -\zeta^{2})/X,~ Y \mapsto \zeta Y/X^{5}.
$$
We also consider an automorphism $\eta'$ of the elliptic curve $E$ defined by
$$
x \mapsto -x-1,~y \mapsto \zeta^{2}y.
$$
Then, we have a morphism $\varphi'$ from $C$ to $E$ which is defined by
$$
\varphi' = \eta' \circ T_{P_{-1}}\circ \varphi\circ \eta : C \longrightarrow E.
$$
The morphism $\varphi'$ is concretely given by
$$
       x = \zeta^{2}X^{3}/(X^{2} - 1),~y= -\zeta^{3}XY/(X^{2} - 1)^{2},
$$
and it is easy to see that
$$
\varphi_{*}'(\frac{\partial}{\partial Y}) = \zeta^{3}X^{3}\frac{\partial}{\partial y} ~\mbox{at any point except the point of infinity.}
$$
Moreover, the homomorphism $\varphi_{*}'$ is injective on the tangent space
at the point at infinity.

Now, we set
$$
\psi = (\varphi, \varphi') : C \longrightarrow A = E_{1} \times E_{2}. 
$$
with $E_{1} = E_{2} = E$. By a direct calculation, we see that 
this morphism is injective. Since $\varphi$ is
\'etale except at the point at infinity and $\varphi'$ is \'etale 
at the point at infinity, we see that the tangent map induced by the morphism 
$\psi$ is injective. Hence, the morphism $\psi$ is an immersion.
We  have
$$
\begin{array}{l}
\varphi(\infty) = P_{\infty}, \\
\varphi (\{(-1, 0), (\zeta^{5}, 0), (\zeta^7, 0)\}) = P_{1},\\
\varphi (\{(1, 0), (\zeta, 0), (\zeta^3, 0)\}) = P_{-1},\\
\varphi (\{(0, 0), (\zeta^{2}, 0), (\zeta^6, 0)\}) = P_{0},
\end{array}
$$
and
$$
\begin{array}{l}
\varphi' (\{\infty, (1, 0), (-1, 0)\}) = P_{\infty}, \\
\varphi' (\{(\zeta, 0), (\zeta^{2}, 0), (\zeta^{7}, 0)\}) = P_{1}\\
\varphi' (\{(\zeta^{3},0), (\zeta^{5}, 0), (\zeta^{6}, 0)\}) = P_{-1},\\
\varphi' ((0, 0)) = P_{0}.
\end{array}
$$
Therefore,
the image of the set $C({\bf F}_{9})$ by the morphism $\psi$ is given by
ten 2-torsion points:
$$
\begin{array}{r}
\{(P_{\infty}, P_{\infty}), (P_{0}, P_{0}), (P_{-1}, P_{\infty}), (P_{-1}, P_{1}),
(P_{0}, P_{1}),\\
 (P_{-1}, P_{-1}), (P_{1}, P_{\infty}), (P_{1}, P_{-1}),
(P_{0}, P_{-1}), (P_{1}, P_{1})\}.
\end{array}
$$
At the ${\bf F}_{9}$-rational points of $C$, the tangent maps are given 
as follows:
$$
\varphi_{*}(\frac{\partial}{\partial Y}) = 0,~
\varphi'_{*}(\frac{\partial}{\partial Y}) \neq 0\quad \mbox{at the point of infinity}
$$
and
$$
\begin{array}{c}
 \varphi_{*}(\frac{\partial}{\partial Y}) = \frac{\partial}{\partial y},~
\varphi_{*}'(\frac{\partial}{\partial Y})= 0\quad \mbox{at the point}~(0, 0)\\
 \varphi_{*}(\frac{\partial}{\partial Y}) = \frac{\partial}{\partial y},~
\varphi_{*}'(\frac{\partial}{\partial Y})= \zeta^{3i + 3}\frac{\partial}{\partial y}\quad \mbox{at the point}~(\zeta^{i}, 0) \\
(i = 0, 1, \ldots, 7).
\end{array}
$$
\begin{lm}
${\rm Im}~\psi$ does not contain any element of 
${\rm Ker} [4]_{A} \backslash {\rm Ker} [2]_{A}$.
Moreover, for any $2$-torsion point $a$, the curve $T_{a}{\rm Im}\psi$ 
does not contain any element of 
${\rm Ker} [4]_{A} \backslash {\rm Ker} [2]_{A}$.
\end{lm}
\proof{
Since the morphism $\psi$ is defined over ${\bf F}_{9}$, the ${\bf F}_{9}$-rational
points on ${\rm Im}~\psi$ come from ${\bf F}_{9}$-rational points on $C$.
Considering ${\rm Ker} [4]_{A} = A({\bf F}_{9})$,
${\rm Ker} [2]_{A} = A({\bf F}_{3})$ and $\psi (C({\bf F}_{9})) \subset A({\bf F}_{3})$,
we see $C$ does not contain elements of 
${\rm Ker} [4]_{A} \backslash {\rm Ker} [2]_{A}$. 

Since 2-torsion points are defined over ${\bf F}_{3}$, 
the latter part follows easily from the former part.}

From here on, we identify $\psi(C)$ with $C$ if we have no confusion. 
We now calculate the intersection number of the curve $C$
with the basis 
$E_{1}, ~E_{2},~ \Delta =\Delta_{1, 1}, ~\Delta_{1, \tau},~\Delta_{1, -\sigma},~\Delta_{1, -\tau\sigma}$. 
Since the degrees of the morphism $\varphi$, $\varphi'$ are equal to 3, we see that $(C, E_{1}) = (C, E_{2}) = 3$. 
Since the genus of $C$ is equal to 4, we have $C^{2} = 6$. 
Calculating the ideals of the pull-backs of the divisors
to $C$, we have the following result.
\begin{lm}
$$
\begin{array}{l}
C^{2} = 6,   \\
(C, E_{1}) = 3, ~(C, E_{2}) = 3,\\
(C, \Delta) =6,~(C, \Delta_{1, \tau}) = 6,\\
(C, \Delta_{1, -\sigma}) = 3,~(C, \Delta_{1, -\tau\sigma}) = 3.
\end{array}
$$
\end{lm}

\begin{th}
In $NS(A)$, we have a decomposition
$$
C = E_{1} + E_{2} - \Delta - \Delta_{1, \tau} + 2\Delta_{1, -\sigma} + 2 \Delta_{1, -\tau\sigma}.
$$
In particular,
$$
j(C) = \left(
\begin{array}{cc}
3  & -(1 + \tau)(1 + 2\sigma) \\
-(1 + 2\sigma^{2})(1 - \tau)  & 3
\end{array}
\right).
$$
\end{th}
\proof{Using the basis $\{E_{1}$, $E_{2}$, $ \Delta$, $\Delta_{1, \tau}$, $\Delta_{1, -\sigma}$,
$\Delta_{1, -\tau\sigma}\}$ of $NS(A)$, we suppose that $C$ is expressed as
$$
C = aE_{1} + bE_{2}+ c\Delta + d\Delta_{1, \tau} + e\Delta_{1, -\sigma} + f\Delta_{1, -\tau\sigma}
$$
with integers $a, b, c, d, e, f$. Considering the intersections of $C$ with the elements of the basis,
we have equations:
$$
\begin{array}{l}
3 = b + c + d + e + f \\
3 = a + c + d + e + f \\
6 = a + b + 2d + e + f \\
6 = a + b + 2c + 2e + f \\
3 = a + b + c + 2d  + 2f \\
3 = a + b + 2c + d + 2e.
\end{array}
$$
Solving these equations, we get the former part of the result.
The latter part follows from Lemma~\ref{basis}.}

\section{Construction of curves on an abelian surface}
We use the notation in previous sections.  Let $E$ be the supersingular elliptic curve
in characteristic 3,
and we consider the superspecial abelian surface
$$
A = E_{1} \times E_{2}
$$
with $E_{1} = E_{2} = E$.

At each point of $\Ker [2]_{A}$, we have 10 tangent directions
which are ${\bf F}_{9}$-rational. At the zero point of $A$ and to 
one of ${\bf F}_{9}$-rational tangent direction at the zero point, 
we construct, in this section,  one curve of genus 4 and 
two elliptic curves which pass through the 
2-torsion point  and  whose tangets at the zero point are equal to 
the tangent direction. At the zero point, we have in total 10 curves of genus 4
and 20 elliptic curves. By translations by 2-torsion points, we have finally 16
curves of genus 4 and 80 elliptic curves. In the next section, we will show
that these curves of genus 4 together with 16 exceptional curves
make a $(16)_{10}$-configuration of rational curves on a Kummer surface,
and the 16 curves of genus 4 and the 80 elliptic curves together with
16 exceptional curves make 112 rational curves on a Kummer surface
which correspond to Leech roots on the Kummer surface.

We  set
$$
     C_{\infty} = C = {\rm Im}~ \psi.
$$
Then, $C_{\infty}$ is a curve of genus 4 and contains the following ten 2-torsion points:
$$
\begin{array}{cl}
 C_{\infty} \ni & (P_{\infty}, P_{\infty}), (P_{1}, P_{\infty}), (P_{-1}, P_{\infty}), (P_{1},  P_{1}), (P_{-1}, P_{1}), (P_{0}, P_{1}), \\
   & (P_{1}, P_{-1}), (P_{-1}, P_{-1}), (P_{0}, P_{-1}), (P_{0}, P_{0}).
\end{array}
$$
We set
$$
\begin{array}{lll}
C_{0} = T_{(P_{0}, P_{0})}^{*}C_{\infty}, &C_{1} = T_{(P_{1}, P_{1})}^{*}C_{\infty}, &C_{\zeta} = T_{(P_{0}, P_{1})}^{*}C_{\infty},\\
C_{\zeta^{2}} = T_{(P_{1}, P_{-1})}^{*}C_{\infty},
&C_{\zeta^{3}} = T_{(P_{-1}, P_{\infty})}^{*}C_{\infty},
&C_{-1} = T_{(P_{-1}, P_{-1})}^{*}C_{\infty},\\
C_{-\zeta} = T_{(P_{0}, P_{-1})}^{*}C_{\infty}, &C_{-\zeta^{2}} = T_{(P_{-1}, P_{1})}^{*}C_{\infty}, &C_{-\zeta^{3}} = T_{(P_{1}, P_{\infty})}^{*}C_{\infty}.
\end{array}
$$
We denote by $T_{A}$ the tangent space of $A$ at the zero point of $A$. 
Note that $T_{A}$ has
a natural basis by using an invariant vector field $\frac{\partial}{\partial y}$
on the elliptic curve $E$.
Then, these 10 curves of genus 4 all pass through the zero point of $A$, and
the tangent vector to the curve $C_{\alpha}$ at the zero point of $A$
gives a point $(1, \alpha)$ in the projectivized tangent space ${\bf P}^{1}(T_{A})$.
They give the ${\bf F}_{9}$-rational points on ${\bf P}^{1}(T_{A})$.
We also note that these 10 curves give the same divisor in $NS(A)$.

Now, to construct 20 elliptic curves which pass through the zero point of $A$,
in the maximal order ${\cal O}={\rm End(E)}$ of $B={\rm End}^{0}(E)$ 
we consider the following 10 decompositions:
$$
\begin{array}{cl}
-(1 + \tau)(1 + 2\sigma)  &  =0 + \bar{\pi}F \\
      & = 0 + \bar{V}\pi   \\
      &   = \overline{(\sigma + \tau)} (-\sigma^2 -\tau)+ \bar{1}\cdot(-\sigma)\\
      &   =\overline{\pi}(-\sigma)  + \bar{1}\cdot\bar{\pi}\sigma^{2}\\
       &  =  \overline{(1 + \sigma^2\tau)}(\tau - \sigma) + \bar{1}\cdot \tau\sigma^2\\
      &   = \bar{1}\cdot\sigma^2\pi + \overline{-{\bar{\pi}}\sigma^2}\cdot1\\
      &   =\overline{(\sigma^{2} +\tau)}(\sigma + \tau)+ \bar{1}\cdot\sigma^2\\
       &  =  \bar{1}\cdot(-\bar{\pi}\sigma) + \overline{\sigma\pi}\cdot 1\\
      &   = \overline{(-\sigma^2 + \tau)}(1 +\tau\sigma^2)+ \bar{1}\cdot(-\tau\sigma) \\
      &   = \overline{(\bar{\pi}\sigma)} \cdot1 + \bar{1}\cdot(-\sigma\pi). \\
\end{array}
$$
Corresponding to these decompositions, we construct the following pairs
of elliptic curves:
$$
\begin{array}{ll}
{\rm class} (\infty) & E_{2} + \Delta_{\pi, F}  \\
{\rm class} (0) & E_{1} + \Delta_{V, \pi}  \\
{\rm class} (1) & \Delta_{\sigma + \tau,-\sigma^2 -\tau} + \Delta_{1, -\sigma} \\
{\rm class} (\zeta) & \Delta_{\pi, -\sigma} + \Delta_{1, \bar{\pi}\sigma^{2}} \\
{\rm class} (\zeta^{2}) & \Delta_{1 + \sigma^2\tau, \tau - \sigma} + 
\Delta_{1, \tau\sigma^2} \\
{\rm class} (\zeta^{3}) &  \Delta_{1, \sigma^2\pi } +
 \Delta_{-{\bar{\pi}\sigma^2}, 1} \\
{\rm class} (-1) &  \Delta_{\sigma^{2} +\tau, \sigma + \tau} + 
\Delta_{1, \sigma^2} \\
{\rm class} (-\zeta) &  \Delta_{1, -\bar{\pi}\sigma} + \Delta_{\sigma\pi, 1} \\
{\rm class} (-\zeta^{2}) &  \Delta_{-\sigma^2 + \tau, 1 +\tau\sigma^2} + 
\Delta_{1,  -\tau\sigma} \\
{\rm class} (-\zeta^{3}) &  \Delta_{\bar{\pi}\sigma, 1} + \Delta_{1, -\sigma\pi} 
\end{array}
$$
In each class ($\alpha$), the tangent vector at the zero point to each elliptic curve
give a point $(1, \alpha)$ in the projectivized tangent space ${\bf P}^{1}(T_{A})$,
and by our construction the divisors which are the sum of two elliptic curves are all
algebraically equivalent to the curve $C_{\alpha}$ of genus 4. We denote by $\cal{E}_{\infty}$
the set of these 20 elliptic curves.

In class ($\alpha$) we denote by $\Delta_{\alpha}$ and  $\Delta_{\alpha}'$
the two ellptic curves in the class. Then, for the triple $\{C_{\alpha}, \Delta_{\alpha}, \Delta_{\alpha}'\}$, we have
$$
C_{\alpha} \equiv  \Delta_{\alpha} + \Delta_{\alpha}'.
$$
By theorem~\ref{intersection}, we can compute the intersection numbers
of these twenty elliptic curves. In particular we see 
$$
     (\Delta_{\alpha}, \Delta_{\alpha}') = 3
$$
in each case. Therefore, we also have 
$$
(C_{\alpha},  \Delta_{\alpha})= 3,\quad  (C_{\alpha},  \Delta_{\alpha}') = 3.
$$  
Since the directions of the tangent vectors of these three curves at the zero point of $A$
coincide with each other, we see that the multiplicity
of the intersection of $\Delta_{\alpha}$ and  $\Delta_{\alpha}'$ at the zero point
is greater than or equal to 2. Since $\Delta_{\alpha}\cap \Delta_{\alpha}'$ 
is a subgroup scheme,
we see $\Delta_{\alpha}\cap \Delta_{\alpha}' = \alpha_{3}$, and so
$\Delta_{\alpha}$ and  $\Delta_{\alpha}'$ intersect with each other at the zero point 
with multiplicity 3, and they don't intersect with each other at any other points.
\begin{lm}
The curve $C_{\alpha}$ intersects wih  $\Delta_{\alpha}$ $($resp. $\Delta_{\alpha}'$$)$
only at the zero point and the multiplicity of their intersection at the
zero point is equal to $3$.
\end{lm}
\proof{Suppose that $C_{\alpha}$ intersects wih  $\Delta_{\alpha}$
at a point $P$ which is not the zero point. Since both $C_{\alpha}$ and $\Delta_{\alpha}$
are defined over ${\bf F}_{9}$, the Galois group over ${\bf F}_{9}$ acts
on rational points on $C_{\alpha}$ and $\Delta_{\alpha}$. In our case
the intersection multiplicity of $C_{\alpha}$ and $\Delta_{\alpha}$ at the zero point
(resp. at the point $P$) is two (resp. one), the Galois group fix the point $P$. 
So $P$ is a ${\bf F}_{9}$-rational point of $C_{\alpha}$. However, $C_{\alpha}$
does not pass ${\bf F}_{9}$-rational points of $A$ which are not ${\bf F}_{3}$-rational points. Therefore, we have $P \in A({\bf F}_{3}) = {\rm Ker}[2]_{A}$. Since we see
that $C_{\alpha}$ does not intersect with $\Delta_{\alpha}$ on ${\rm Ker}[2]_{A}\setminus \{0\}$, we conclude that $C_{\alpha}$ intersects wih  $\Delta_{\alpha}$ only at the zero point with multiplicity 3. The rest is similarily proved.}

Take a 2-torsion point $a$ which is not the zero point. Then, we can check that
just one curve among 3 curves $C_{\alpha}$, $\Delta_{\alpha}$ and $\Delta_{\alpha}'$
passes through the 2-torsion point $a$.

Now, let $\{C_{\beta}, \Delta_{\beta}, \Delta_{\beta}'\}$ be a  triple different
from $\{C_{\alpha}, \Delta_{\alpha}, \Delta_{\alpha}'\}$.
Then, we have
$$
     3 =  ( C_{\alpha},   \Delta_{\beta})
 = (\Delta_{\alpha}, \Delta_{\beta}) + (\Delta_{\alpha}', \Delta_{\beta}).
$$
Therefore, we have $ (\Delta_{\alpha}, \Delta_{\beta}) = 0$, $1$, $2$ or $3$.
If $ (\Delta_{\alpha}, \Delta_{\beta}) = 0$, then $\Delta_{\beta}$ is a translation of $\Delta_{\alpha}$. However, since both divisors contain the zero point, we see
$\Delta_{\beta}=\Delta_{\alpha}$, which contradicts to $\alpha \neq \beta$.
If $ (\Delta_{\alpha}, \Delta_{\beta}) = 3$, we have 
$ (\Delta_{\alpha}', \Delta_{\beta}) = 0$. Therefore, by the same argument
we have $\Delta_{\beta}=\Delta_{\alpha}'$, a contradiction. Hence we have
the following.
\begin{lm}
$ (\Delta_{\alpha}, \Delta_{\beta})$ is equal to either $1$ or $2$ 
if $\alpha \neq \beta$.
\end{lm}
In case $ (\Delta_{\alpha}, \Delta_{\beta}) = 1$, $\Delta_{\alpha}$ 
intersects with $\Delta_{\beta}$ only at the zero point transversely.
In case $ (\Delta_{\alpha}, \Delta_{\beta}) = 2$, 
since  $\Delta_{\beta}\cap \Delta_{\alpha}$ is a subgroup scheme of $A$, $\Delta_{\alpha}$ intersects with $\Delta_{\beta}$ at the zero point  
and a point of order two transversely.

For the group $\Ker [2]_{A}$ of 2-torsion points, we set
$$
{\cal E} = \{T_{a}G~\mid ~G\in {\cal E}_{\infty},  a \in \Ker [2]_{A}\}.
$$
Since each elliptic curve $G$ contains four 2-torsion points, we see that
the set ${\cal E}$ contains $20 \times 16\div 4 = 80$ elliptic curves.
We denote by ${\cal E}_{a}$ the subset of ${\cal E}$
such that $G \in {\cal E}$ contains the point $a \in \Ker [2]_{A}$.
We have ${\cal E}_{0} = {\cal E}_{\infty}$.
Then, the set ${\cal E}_{a}$ contains 20 elliptic curves whose tangent vectors 
at the point $a$ give ten ${\bf F}_{9}$-rational points of the projectivized
tangent space at the point $a$. Two of them give the same rational point
in the projectivised tangent space.

We set 
$$
\begin{array}{l}
Q_{1} =  (\zeta, \zeta^{3}),~Q_{2} =  (\zeta^{2}, \zeta),
~Q_{3} =  (\zeta^{3}, \zeta),\\
Q_{4} =  (-\zeta, \zeta),~Q_{5} =  -(\zeta^{2}, \zeta^{3}),
~Q_{6} =  (-\zeta^{3}, \zeta^{3}).
\end{array}
$$
Then, the 4-torsion points on $E$ which are not in the group $\Ker [2]_{E}$
of 2-torsion points are given by
$$
Q_{1}, Q_{2}, Q_{3}, Q_{4}, Q_{5}, Q_{6}, -Q_{1}, 
-Q_{2}, -Q_{3}, -Q_{4}, -Q_{5}, -Q_{6}.
$$
and the 2-torsion points are given by
$$
P_{\infty}, P_{0}, P_{1}, P_{-1}.
$$
For each elliptic curve  in ${\cal E}_{0}$ we put a number  as follows:
$$
\begin{array}{ll}
(1)~\Delta_{V, \pi} , ~                           (1)'~  E_{1}, &
(2)~  E_{2} , ~                                       (2)' ~ \Delta_{\pi, F},\\
(3)~  \Delta_{\sigma^{2} \tau, \sigma + \tau},~     (3)' ~\Delta_{1, \sigma^2} &
(4)~ \Delta_{\pi, -\sigma} ,~                      (4)' ~   \Delta_{1, \bar{\pi}\sigma^{2}},\\
(5)~ \Delta_{-\sigma^2 + \tau, 1 +\tau\sigma^2},~    (5)' ~  \Delta_{1,  -\tau\sigma}, &
(6)~ \Delta_{1, \sigma^2\pi },~                   (6)' ~   \Delta_{-{\bar{\pi}\sigma^2}, 1},\\
(7)~ \Delta_{\sigma + \tau,-\sigma^2 -\tau},~    (7)' ~ \Delta_{1, -\sigma},&
(8) ~  \Delta_{1, -\bar{\pi}\sigma} ,~             (8)' ~  \Delta_{\sigma\pi, 1} ,\\
(9)~  \Delta_{1 + \sigma^2\tau, \tau - \sigma},~     (9)' ~ \Delta_{1, \tau\sigma^2} ,&
(10)~  \Delta_{\bar{\pi}\sigma, 1},~              (10)' ~   \Delta_{1, -\sigma\pi} .
\end{array}
$$
Then, we have the following list of elliptic curves in ${\cal E}_{0}$ which pass through
the 4-torsion points on $A = E_{1}\times E_{2}$ which are not in the group 
$\Ker [2]_{A}$. Since the table is very big, we divide it into two parts.
Here, the points in the top row of the table  are the ones on $E_{1}$, and
the points in the left-end column are the ones in $E_{2}$. For example
this table shows that the 4-torsion point $(Q_{1}, Q_{1}) \in A$ lies on the elliptic curve
$(9)~ \Delta_{1 + \sigma^2\tau, \tau - \sigma}$, and
the 4-torsion point $(Q_{2}, Q_{1}) \in A$
lies on the elliptic curve $(9)' ~ \Delta_{1, \tau\sigma^2}$.
From this table we see that each 4-torsion point of $A$ which is not in $\Ker [2]_{A}$
is on only one elliptic curve in ${\cal E}_{0}$.

\begin{center}
\begin{tabular}{| l |c|c|c|c|c|c|c|c|}
\hline
     & $Q_{1} $ & $Q_{2}$ & $Q_{3}$ & $Q_{4}$ & $Q_{5}$ & $Q_{6}$ & $-Q_{1}$ & $-Q_{2}$ \\
\hline
$Q_{1}$     & $(9) $ & $(9)'$ & $(10)'$ & $(8)'$ & $(1)$ & $(7)'$ & $(9)$ & $(8)$ \\
\hline
$Q_{2}$     & $(10) $ & $(6)'$ & $(7)'$ & $(2)'$ & $(8)$ & $(5)'$ & $(9)'$ & $(4)'$ \\
\hline
$Q_{3}$     & $(5)' $ & $(8)$ & $(10)'$ & $(7)'$ & $(1)$ & $(3)$ & $(8)'$ & $(3)'$ \\
\hline
$Q_{4}$     & $(4) $ & $(7)'$ & $(4)$ & $(5)$ & $(5)'$ & $(6)$ & $(10)'$ & $(1)$ \\
\hline
 $Q_{5}$    & $(7)' $ & $(4)'$ & $(9)'$ & $(6)'$ & $(8)$ & $(6)'$ & $(2)'$ & $(10)$ \\
\hline
 $Q_{6}$    & $(10)' $ & $(1)$ & $(7)$ & $(9)'$ & $(7)'$ & $(4)$ & $(3)'$ & $(5)'$ \\
\hline
$-Q_{1}$     & $(9) $ & $(8)$ & $(5)'$ & $(6)$ & $(3)'$ & $(8)'$ & $(9)$ & $(9)'$ \\
\hline
$-Q_{2}$     & $(9)'$ & $(4)'$ & $(10)$ & $(3)'$ & $(6)'$ & $(2)'$ & $(10)$ & $(6)'$ \\
\hline
$-Q_{3}$     & $(8)' $ & $(3)'$ & $(8)'$ & $(6)$ & $(9)'$ & $(3)$ & $(5)'$ & $(8)$ \\
\hline
$-Q_{4}$     & $ (10)'$ & $(1)$ & $(3)'$ & $(5)$ & $(4)'$ & $(9)'$ & $(4)$ & $(7)'$ \\
\hline
$-Q_{5}$     & $(2)' $ & $(10)$ & $(2)'$ & $(5)'$ & $(10)$ & $(3)'$ & $(7)'$ & $(4)'$ \\
\hline
$-Q_{6}$     & $(3)' $ & $(5)'$ & $(7)$ & $(4)$ & $(4)'$ & $(6)$ & $(10)'$ & $(1)$ \\
\hline
$P_{\infty}$     & $(1)' $ & $(1)' $ & $(1)' $ & $(1)' $ & $(1)' $ & $(1)' $ & $(1)' $ & $(1)' $ \\
\hline
$P_{0}$     & $(7) $ & $(8)'$ & $(9)$ & $(3)$ & $(4)$ & $(5)$ & $(7)$ & $(8)'$ \\
\hline
$P_{1}$     & $(7) $ & $(4)$ & $(6)'$ & $(10)$ & $(2)'$ & $(5)$ & $(7)$ & $(4)$ \\
\hline
$P_{-1}$     & $(6)' $ & $(2)'$ & $(9)$ & $(3)$ & $(8)'$ & $(10)$ & $(6)'$ & $(2)'$ \\
\hline
\end{tabular}
\end{center}

\begin{center}
\begin{tabular}{| l |c|c|c|c|c|c|c|c|}
\hline
    & $ -Q_{3}$ & $-Q_{4}$ & $-Q_{5}$ & $-Q_{6}$ & $P_{\infty}$ & $P_{0}$ & $P_{1}$ & $P_{-1}$ \\
\hline
$Q_{1}$     & $(5)'$ & $(6)$ & $(3)'$ & $(8)'$ & $(2)$ & $(3)$ & $(3)$ & $(4)'$ \\
\hline
$Q_{2}$     & $(10) $ & $(3)'$ & $(6)'$ & $(2)'$ & $(2)$ & $(10)'$ & $(6)$ & $(1)$ \\
\hline
$Q_{3}$     & $(8)' $ & $(6)$ & $(9)'$ & $(3)$ & $(2)$ & $(9)$ & $(4)'$ & $(9)$ \\
\hline
$Q_{4}$     & $(3)'$ & $(5)$ & $(4)'$ & $(9)'$ & $(2)$ & $(7)$ & $(8)$ & $(7)$ \\
\hline
 $Q_{5}$    & $(2)' $ & $(5)'$ & $(10)$ & $(3)'$ & $(2)$ & $(6)$ & $(1)$ & $(10)'$ \\
\hline
$Q_{6}$     & $(7) $ & $(4)$ & $(4)'$ & $(6)$ & $(2)$ & $(5)$ & $(5)$ & $(8)$ \\
\hline
$-Q_{1}$     & $(10)' $ & $(8)'$ & $(1)$ & $(7)'$ & $(2)$ & $(3)$ & $(3)$ & $(4)'$ \\
\hline
$-Q_{2}$     & $(7)'$ & $(2)'$ & $(8)$ & $(5)'$ & $(2)$ & $(10)'$ & $(6)$ & $(1)$ \\
\hline
$-Q_{3}$     & $(10)'$ & $(7)'$ & $(1)$ & $(3)$ & $(2)$ & $(9)$ & $(4)'$ & $(9)$ \\
\hline
$-Q_{4}$     & $(4)$ & $(5)$ & $(5)'$ & $(6)$ & $(2)$ & $(7)$ & $(8)$ & $(7)$ \\
\hline
$-Q_{5}$     & $(9)'$ & $(6)'$ & $(8)$ & $(6)'$ & $(2)$ & $(6)$ & $(1)$ & $(10)'$ \\
\hline
$-Q_{6}$     & $(7) $ & $(9)'$ & $(7)'$ & $(4)$ & $(2)$ & $(5)$ & $(5)$ & $(8)$ \\
\hline
$P_{\infty}$     & $(1)' $ & $(1)'$ & $(1)'$ & $(1)'$ & \multicolumn{4}{|c|}{} \\
\cline{1-5}
$P_{0}$     & $(9) $ & $(3)$ & $(4)$ & $(5)$ & \multicolumn{4}{|c|}{2-torsion points} \\
\cline{1-5}
$P_{1}$     & $ (6)'$ & $(10)$ & $(2)'$ & $(5)$ & \multicolumn{4}{|c|}{} \\
\cline{1-5}
$P_{-1}$     & $(9) $ & $(3)$ & $(8)'$ & $(10)$ & \multicolumn{4}{|c|}{} \\
\hline
\end{tabular}
\end{center}

As for the elliptic curves in ${\cal E}$, any curve in ${\cal E}$ is a translation
of an elliptic curve in ${\cal E}_{0}$ by a 2-torsion point.
Since  there exist sixteen points in $\Ker [2]_{A}$ and each elliptic curve
contains four 2-torsion points, we see that
each 4-torsion point of $A$ which is not in $\Ker [2]_{A}$
is on just 4 elliptic curves in ${\cal E}$.   Each two of these 4 elliptic curves intersects 
transversely at the 4-torsion point with each other. 
The inversion of $A$ acts on sixteen 4-torsion points of an elliptic curve $G$ 
in ${\cal E}$, and it fixes four 2-torsion points. Therefore, it makes 10 orbits.
We will see later that these 10 orbits make ten ${\bf F}_{9}$-rational points 
on the rational curve on the Kummer surface ${\rm Km}(A)$ which
is obtained as the image of the elliptic curve $G$.

We set now
$$
\begin{array}{r}
S = \{(P_{\infty}, P_{\infty}), (P_{1}, P_{\infty}), (P_{-1}, P_{\infty}), (P_{1}, P_{1}), (P_{-1}, P_{1}), (P_{0}, P_{1}), \\
\quad \quad (P_{1}, P_{-1}), (P_{-1}, P_{-1}),
(P_{0}, P_{-1}), (P_{0}, P_{0})\}
\end{array}
$$
and  set
$$
      {\cal C} = \{T_{a}C_{\infty} \mid a \in S\}.
$$
Then, we have ten curves of genus 4 which pass through the origin of $A$.
As we already stated, the tangent vectors to these curves in ${\cal C}$ 
at the origin give
ten ${\bf F}_{9}$-rational points of the projectivized tangent space ${\bf P}^{1}$
at the origin of $A$.
More generally, for the group $\Ker [2]_{A}$ of 2-torsion points, we set
$$
{\cal D} = \{T_{a}C_{\infty}\mid a \in \Ker [2]_{A}\}.
$$
Then, we have sixteen curves of genus 4 such that each passes through
10 points in $\Ker [2]_{A}$. 
Since $(T_{a}C_{\infty}\cdot T_{b}C_{\infty}) = 6$ for any 2-torsion points $a, b$, 
the curves $T_{a}C_{\infty}$ and $T_{b}C_{\infty}$ intersect each other 
at six 2-torsion points transversely and they don't intersect at any other points.
Therefore, two of the elements in ${\cal D}$ intersect at 6 points in $\Ker [2]_{A}$
transversely. We set
$$
{\cal C}_{b} = \{T_{b + a}C_{\infty} \mid a \in S\} \subset {\cal D}.
$$
Then, ${\cal C}_{b}$ is a subset of ${\cal D}$ and  we have 10 curves 
of genus 4 which pass through the point $b$. Note that
the other six curves in ${\cal D}$ don't pass the point $b$.
Thus, we get the set ${\cal D}$ of 16 curves of genus 4 and the set ${\cal E}$
of 80 elliptic curves. Summarizing these results, we have the following theorem.
\begin{th}
Under the notation above, on the superspecial abelian surface $A$, 
we have the set ${\cal D}$ of $16$ curves of genus $4$ and 
the set ${\cal E}$ of $80$ elliptic curves. 

{\rm (i)} Different two curves $C_{1}$, $C_{2}$ in ${\cal D}$ intersect each other 
only at $6$ points in $\Ker [2]_{A}$. Any curve in ${\cal D}$ doesn't pass
through any point in $\Ker [4]_{A}\setminus \Ker [2]_{A}$.

{\rm (ii)} Different two curves $C_{1}$, $C_{2}$ in ${\cal E}$ don't intersect 
each other except at points in $\Ker [4]_{A}$.

{\rm (iii)} A curve $C_{1}$ in ${\cal D}$ and a curve $C_{2}$ in ${\cal E}$ don't
 intersect each other except at points in $\Ker [2]_{A}$.

{\rm (iv)} At a point $a$ in $\Ker [2]_{A}$
and to a ${\bf F}_{9}$-rational tangent direction $t$ at the point $a$,
we have a triple of a curve $C$ in ${\cal D}$ and two elliptic curves $\Delta$,
 $\Delta'$ in ${\cal E}$
such that these three curves pass through the point $a$ and whose tangent
directions at the point $a$ coincide with the direction $t$. These three curves
intersect each other only at the point $a$, and the intersection numbers
of two curves among these three are all equal to $3$. Moreover, 
$C$ is algebraically equivalent to $\Delta + \Delta'$.

{\rm (v)} For a point $b$ in
$\Ker [4]_{A}\setminus \Ker [2]_{A}$, there exist just $4$ elliptic curves in 
${\cal E}$ which pass through the point $b$. If two elliptic curves
in ${\cal E}$ intersect each other at the point $b$, then they also
intersect each other at the point $-b$. These $4$ elliptic curves
intersect each other transversely at the points $b$ and $-b$. 
\end{th}

\section{Rational curves on a Kummer surface}\label{kummer}
Let $A = E_{1} \times E_{2}$ with $E_{1} = E_{2} = E$ 
be the superspecial abelian surface and let $\iota$ be the inversion of $E$. 
We consider the quotient surface
$A/\langle \iota \times \iota \rangle$ and its minimal resolution 
${\rm Km}(A)$. ${\rm Km}(A)$ is a K3 surface and is called a Kummer surface.
We consider the following diagram:
$$
\begin{array}{ccc}
{\tilde A}& \stackrel{f}{\longrightarrow} &  A\\ 
  \quad  \downarrow{g}  &  &        \downarrow\\  
{\rm Km}(A) & \stackrel{f'}{\longrightarrow} &   A/\langle \iota \times \iota \rangle .   
\end{array}
$$
Here, $f'$ is the minimal resolution of singularities and $f$ is the blowings-up
at sixteen 2-torsion points of $A$. The morphism $g$ is the quotient map
to the quotient surface by the group of order 2.

Let $\{C_{\alpha}, \Delta_{\alpha}, \Delta_{\alpha}'\}$ be a triple of
the curve of genus 4 and two elliptic curves $\Delta_{\alpha}, \Delta_{\alpha}'$
as in Section 8 which pass through the zero point of $A$. We denote by $\ell_{0}$
the exceptional curve which comes from the zero point of $A$. We denote by
$\tilde{C}_{\alpha}$, $\tilde{\Delta}_{\alpha}$, $\tilde{ \Delta}_{\alpha}$
the proper transforms of $C_{\alpha}$, $\Delta_{\alpha}$, $\Delta_{\alpha}$
respectively. Then, each two of members of the tirple intersects
at the origin with multiplicity 3, we see that each two of the proper transforms
intersects on the exceptinal curve $\ell_{0}$ with multiplicity 2.
Since $\Delta_{\alpha}$ and $ \Delta_{\alpha}'$ 
contain four points  in $\Ker [2]_{A}$ respectively
and $C_{\alpha}$ contains 10 points in $\Ker [2]_{A}$,
we have the following intersection numbers for the proper transforms:
$$
\begin{array}{l}
\tilde{C}_{\alpha}^{2} =  -4,~  \tilde{\Delta}_{\alpha}^2 = -4, ~
\Delta_{\alpha}'^{2} = -4,\\
(\tilde{C}_{\alpha},  \tilde{\Delta}_{\alpha}^{2}) = 2,~
(\tilde{C}_{\alpha},  \tilde{\Delta}_{\alpha}'^{2}) = 2, ~
 (\tilde{\Delta}_{\alpha},  \tilde{\Delta}_{\alpha}') = 2\\
\ell_{0}^{2} = -1,~ (\tilde{C}_{\alpha}, \ell_{0}) = 1,~
 (\tilde{\Delta}_{\alpha}, \ell_{0}) = 1,~
(\tilde{\Delta}_{\alpha}', \ell_{0}) = 1.
\end{array}
$$
It is clear that $g(\tilde{C}_{\alpha})$, $g(\tilde{\Delta}_{\alpha})$,
$g(\tilde{ \Delta}_{\alpha})$ and $g(\ell_{0})$ are all nonsingular rational curves.
By the construction we have
$$
g^{-1}(g(\tilde{C}_{\alpha})) = \tilde{C}_{\alpha},~ 
g^{-1}(g(\tilde{\Delta}_{\alpha}))=\tilde{\Delta}_{\alpha}, ~ 
g^{-1}(g(\tilde{ \Delta}_{\alpha}))= \tilde{\Delta}_{\alpha}', ~
g^{-1}(g(\ell_{0}))= 2\ell_{0}.
$$
Hence, by $\deg g = 2$ we conclude
$$
\begin{array}{l}
g(\tilde{C}_{\alpha})^{2} =  -2, ~ g(\tilde{\Delta}_{\alpha})^2 = -2, ~
g(\Delta_{\alpha}')^{2} = -2, ~g(\ell_{0})^{2} = -2\\
(g(\tilde{C}_{\alpha}),  g(\tilde{\Delta}_{\alpha})^{2}) = 1,~
(g(\tilde{C}_{\alpha}),  g(\tilde{\Delta}_{\alpha})'^{2}) = 1,~
 (g(\tilde{\Delta}_{\alpha}),  g(\tilde{\Delta}_{\alpha}')) = 1,\\
(g(\tilde{C}_{\alpha}), g(\ell_{0})) = 1,~ (q(\tilde{\Delta}_{\alpha}), g(\ell_{0})) = 1,~
(g(\tilde{\Delta}_{\alpha}'), g(\ell_{0})) = 1.
\end{array}
$$
We note that the rational curves $g(\tilde{C}_{\alpha})$,  $g(\tilde{\Delta}_{\alpha})$,
$g(\tilde{\Delta}_{\alpha})'$ and $g(\ell_{0})$ intersect one another
at the same point.

We denote by ${\cal A}$ the set of exceptional curves for the resolution $f'$
and by ${\tilde{\cal D}}$ the set of of proper transforms  of the curves 
in ${\cal D}$.
We also denote by ${\cal B}$ the set of rational curves on ${\rm Km(A)}$ 
which are obtained
as the images by $g$ of the proper transforms. Then, from the argument above,
we have the following theorem (see Proposition \ref{16-10}).

\begin{th}\label{16-10th} The sets ${\cal A}$ and ${\cal B}$ make a $(16)_{10}$-configuration
on a supersingular $K3$ surface ${\rm Km}(A)$,
that is, the sets ${\cal A}$ and ${\cal B}$ respectively consist of $16$ non-singular rational 
curves and each curve in one set  intersects exactly $10$ curves from the other set
transversely.
\end{th}
We consider the set ${\tilde{\cal E}}$ of proper transforms  of the curves in ${\cal E}$, 
and
we denote by ${\bar{\cal E}}$ the set of rational curves which are obtained
as the images of the curves in ${\tilde{\cal E}}$ by $g$.
We set 
$$
{\cal R} = {\cal A}\cup {\cal B} \cup {\bar{\cal E}}.
$$
Then, ${\cal R}$ contains in total 112 nonsigular rational curves
whose self-intersection numbers are all equal to $-2$.
By the argument above,
we know that on one rational curve  $\ell$ in ${\cal R}$ there exists ten
${\bf F}_{9}$-rational points and that at each ${\bf F}_{9}$-rational point
4 rational curves in ${\cal R}$ intersect transversely each other.
Therefore we have $10\times 112\div 4 = 280$ ${\bf F}_{9}$-rational points.
We denote by ${\cal P}$ the set of these 280 points.
Summarizing these resuts, we have the following theorem 
(see Proposition \ref{280-4}).
\begin{th}\label{280-4th}
The set ${\cal R}$ is the set of Leech roots in $NS({\rm Km}(A))$. 
Moreover, ${\cal P}$ and 
${\cal R}$ make a $(280_{4}, 112_{10})$-configuration.
\end{th}

\vspace{0.5cm}
\noindent
T.\ Katsura: Faculty of Science and Engineering, Hosei University,
Koganei-shi, Tokyo 184-8584, Japan

\noindent
E-mail address: toshiyuki.katsura.tk@hosei.ac.jp

\vspace{0.3cm}
\noindent
S.\ Kond\=o: Graduate School of Mathematics, Nagoya University, Nagoya 464-8602, Japan

\noindent
E-mail address: kondo@math.nagoya-u.ac.jp

\end{document}